\documentclass[twocolumn]{article}

\usepackage[utf8]{inputenc}
\usepackage[colorinlistoftodos]{todonotes}
\usepackage{xcolor}
\usepackage{lipsum}
\usepackage{amsthm}
\usepackage{amssymb}
\usepackage{amsmath}
\usepackage{endnotes}
\usepackage{enumitem}
\usepackage{mathpazo}

\usepackage{hyperref}
\hypersetup{
    colorlinks= true,
    linkcolor = blue,
    filecolor = blue,
    citecolor = blue,      
    urlcolor  = blue,
    }
\usepackage{cleveref}
\usepackage{url}
\usepackage[backend=biber,defernumbers=true,style=authoryear]{biblatex} 

\usepackage[margin=1in]{geometry}

\newtheorem{gp}{Principle}
\crefname{gp}{principle}{principles}
\Crefname{gp}{Principle}{Principles}

\theoremstyle{definition}

\theoremstyle{definition}
\newtheorem*{dfn*}{Definition}

\newcommand{\affil}{\thanks}

\makeatletter
\def\enoteformat{  \rightskip\z@ \leftskip=3em \parindent=1.8em
  \leavevmode{\setbox\z@=\lastbox}\llap{P\theenmark.\enskip}}
\makeatother

\newcommand\blfootnote[1]{  \begingroup
  \renewcommand\thefootnote{}\footnote{#1}  \addtocounter{footnote}{-1}  \endgroup
}

\DeclareBibliographyCategory{fullcited}

\newcommand{\notescitet}[1]{\textcite{#1}\addtocategory{fullcited}{#1}}

\title{Data Science and Social Justice in the Mathematics Community}
\author{
        Quindel Jones\affil{Quindel Jones is a PhD candidate in systems modeling and analysis at Virginia Commonwealth University.},
        Andr\'es R. Vindas Mel\'endez\affil{Andr\'es R. Vindas Mel\'endez is an NSF postdoctoral scholar at UC Berkeley. In July 2024, he will start as an assistant professor of mathematics at Harvey Mudd College.},
        Ariana Mendible\affil{Ariana Mendible is an assistant professor of mathematics at Seattle University.},
        Manuchehr Aminian\affil{Manuchehr Aminian is an assistant professor of mathematics at California State Polytechnic University, Pomona.},\\ 
        Heather Z. Brooks\affil{Heather Zinn Brooks is an assistant professor of mathematics at Harvey Mudd College.},
        Nathan Alexander\affil{Nathan Alexander is an assistant professor of mathematics and data science at Morehouse College.},
        Carrie Diaz Eaton\affil{Carrie Diaz Eaton is an associate professor of digital and computational studies at Bates College.},
        and Philip Chodrow\affil{Philip Chodrow is an assistant professor of computer science at Middlebury College. His email address is \texttt{pchodrow@middlebury.edu}.}
        }

\date{March 2023}

\begin{document}

\maketitle

\section{Introduction} \label{sec:intro}
\addtoendnotes{\vspace{2em}\noindent {\large \sc Notes on \Cref{sec:intro}}\vspace{1em}}

\blfootnote{The first seven authors are ordered by career stage at the time of article submission.}

In the summer of 2021, the Puerto Rico Association of Criminal Defense Lawyers (PRACDL) submitted a brief to the US Supreme Court in the case \emph{Rodr\'iguez-Rivera vs.~USA}. 
At issue was whether a certain class of drug conspiracy crimes should be considered ``controlled substance offenses,'' a status which can carry enhanced sentences. 
Drug conspiracy charges implicate multiple defendants and carry a low burden of evidence required to obtain a conviction against alleged co-conspirators. 
PRACDL's brief argued that prosecution in drug conspiracy cases reflected discrimination and disparate impact along axes of race, ethnicity, and class. 
The brief was supported by analysis from the Institute for the Quantitative Study of Inclusion, Diversity, and Equity (QSIDE), including work by present authors Manuchehr Aminian and Phil Chodrow. \endnotetext{QSIDE's website is \url{https://qsideinstitute.org/}. }
Working in close dialogue with lawyers at PRACDL, the QSIDE team wrote parsers to extract information from unstructured docket data, combined these results with geographic Census data, and tested hypotheses that conspiracy charges disproportionately took place in areas of Puerto Rico impacted by poverty, unemployment, and low rates of education.
      
This work with PRACDL is a small entry in a growing body of scholarship in \emph{data science for social justice} (DS4SJ).\footnote{We also use the phrase \emph{justice data science} interchangeably to refer to the same scholarly area.}\endnotetext{See \emph{Flores Gonz\'alez vs USA , 2022} for another legal brief by PRACDL supported by data analysis.} 
This area has received growing attention in several overlapping professional communities, including mathematics, statistics, and computer science. 
Our focus in this article is how members of the mathematical community can learn, practice, and support this burgeoning area of scholarship and activism.

We highlight ways that DS4SJ work can both appeal to and challenge trained mathematicians.
Scholars who pursue DS4SJ may look forward to positive outcomes for served communities, new research questions, new professional connections outside the academy, and opportunities to connect their identity or values to their scholarly work. 
Realizing these benefits, however, requires confronting a range of challenges and risks. 
First, impact is never guaranteed. 
In the PRACDL case, for example, the petition to the Court was ultimately denied without comment.  
Furthermore, the process of justice data science calls on skills in which many of us are untrained. 
Realizing impact beyond the ivory tower requires engagement beyond the ivory tower, with affected populations, domain experts, and decision-makers. 
When this work is pursued without sufficient attention to the needs of impacted communities, the results can fail as both scholarship and activism.
Attentive engagement with non-scholarly collaborators is a skill not taught in mathematics training programs, perhaps because it is not viewed as ``mathematical.''
Finally, researchers who take on justice data science work must naturally keep an eye on how it aligns with their needs for advancement in their institutions and professional communities. 

Our overall aim is to help scholars who wish to begin or deepen their engagement with DS4SJ work, and advocate to the broader mathematics community that this work should be valued and supported. 
Our core arguments are: 
\begin{itemize}
    \item[(a).] Members of the mathematical community, broadly construed, have contributions to make towards social justice causes through both scholarship and pedagogy. 
    \item[(b).] Making these contributions requires us to engage with collaborators and partners outside the academy, including impacted communities, policy makers, and others who contribute experience, knowledge, and power beyond our own. 
    \item[(c).] Data science work supporting social justice can indeed be a scholarly, scientific endeavor, which should be supported by mathematical institutions at all levels. 
\end{itemize}
    
In \Cref{sec:DS4SJ} we define data science for social justice (DS4SJ), construing this work broadly to include a wide range of computational, mathematical, and qualitative work. 
We compare and contrast DS4SJ to related bodies of work, including data science for social good and academic data science of human behavior. 
Based on our definition, we propose a set of four guiding principles for the practice of DS4SJ.
These principles describe work that actively engages community members and decision-makers; is aware of its own limitations; intentionally aims to divest from privilege; and grapples with the systemic nature of oppression. 
We continue in \Cref{sec:challenges-opportunities} with discussion of some features of recent and ongoing work in DS4SJ. 
We especially focus on considerations for early-career scholars considering greater involvement this work in either research or pedagogy. 
We also highlight benefits for departments that make a commitment to supporting this work.
In \Cref{sec:vignettes}, we illustrate the practice of DS4SJ with a series of vignettes in which early-career scholars describe their work in DS4SJ and its impact on their growth in the broad mathematical community. 
We conclude in \Cref{sec:conclusion} with a reflection and discussion on the outlook for future work. 

\section{Data Science For Social Justice}\label{sec:DS4SJ}
\addtoendnotes{\vspace{2em}\noindent {\large \sc Notes on \Cref{sec:DS4SJ}}\vspace{1em}}

We take an intentionally broad view of \emph{data science}, referring to the wide set of practices, methods, and activities that support the use of data to gain insight and guide action. 
This includes standard steps in the data analysis pipeline, such as data collection, exploration, visualization, modeling, and communication. 
We also deliberately include an expansive range of other activities, including: formulation of mathematical models, development of data analysis and learning algorithms; qualitative research that contextualizes analysis; and deployment of publicly-facing data products such as visualizations and dashboards. 
In this broad view, data science work takes in many academic disciplines. 
It is not restricted to departments of statistics, mathematics, or computer science; it is indeed not restricted to STEM fields at all. 
That said, while our view of data science is broad, our focus in this article is primarily on how the mathematical community can engage with data science work. 
For this reason, our discussion largely centers on activities and practices in which mathematical training is relevant. 

We frame data science for social justice within our broad construal of data science as a whole. 
\begin{dfn*} \label{def:dssj}
    \emph{Data science for social justice} (DS4SJ) is data scientific work (broadly construed) that actively challenges systems of inequity and concretely supports the liberation of oppressed and marginalized communities.\footnote{Forms of oppression and marginalization include but are not limited to racism, colonialism, classism, sexism, homophobia, transphobia, and ableism.} 
\end{dfn*}
DS4SJ acknowledges that some groups of people have been historically oppressed, marginalized, and disenfranchised, and aims to rectify these harms. 
It does not necessarily aspire to benefit ``society as a whole,'' but specifically aims to correct injustice.
DS4SJ is therefore distinct from data science for social \emph{good}, a phrase that often refers to any data scientific work with net positive social impact.
Work featured by the organization named Data Science For Social Good, for example, includes efforts towards sustainability, efficient public infrastructure, and reduction of political corruption.\endnotetext{Data Science For Social Good's website is \url{https://www.datascienceforsocialgood.org}.}  
These are general social goods which benefit large segments of the world's population. 
Social \emph{justice} applications also benefit populations, but are distinctive in their  explicit focus on the needs of oppressed and marginalized communities harmed by historical wrongs. 
We emphasize that this distinction, though important, is also porous. 
Many social good projects also have social justice elements; for example, the improvement of public infrastructure may have an especially large impact on communities that have been historically underserved by that infrastructure. 
Similarly, justice efforts that rectify systematic wrongs may be appropriately viewed as promoting broader social good. 
While it is important to distinguish these aims and scopes, we emphasize that neither precludes the other.
\endnotetext{We acknowledge that the question of which populations are oppressed or marginalized, and in which ways, may be unclear or contentious. 
Because of this, there may not be consensus on whether a given organizing effort is best framed as aiming towards social justice, social good more broadly, or something else entirely. 
These same considerations apply to data science efforts, and imply that the line between data science for social justice, data science for social good, and other forms of socially-motivated data work may be porous and contested.}

DS4SJ is also importantly related to work on bias and harm in large-scale machine learning systems. 
Recent work has shown that many of these models contain significant biases against members of marginalized identity groups. 
These effects are often intersectional: in a famous study by \textcite{pmlr-v81-buolamwini18a}, facial recognition algorithms performed worse on Black women than they did on either Black men or white women, with the faces of white men being recognized most reliably. 
Contributors to bias in machine learning models include underrepresentation of marginalized identities or overrepresentation of oppressive stereotypes in training data; use of flawed proxy variables for measuring predictive accuracy; or failure to consider the role of an automated system in social and historical context. 
With the growing prevalence of machine learning models in personal technology, hiring, policing, and healthcare, the stakes for designing nondiscriminatory systems at mass scales are extremely high. 
We view work toward fair machine learning as DS4SJ: it is data science that challenges systems of inequity and supports the liberation of oppressed and marginalized communities. 
This is especially true when the work influences the design and assessment of deployed, large-scale systems.
That said, the majority of this work is carried out in the computer science community rather than the mathematics community and is therefore outside our scope. 
We therefore refer the interested reader to \textcite{mehrabi2021survey} for a survey of recent work in bias and fairness in machine learning, including conceptual foundations, empirical findings, and proposed methods for reducing algorithmic discrimination in a range of machine learning tasks. 

Finally, we distinguish DS4SJ as a specific subset of the broad academic area of data science and modeling of social systems. 
This large and growing area includes work in applied mathematics, computational social science, physics, network science, computer science, and allied fields. 
Some of this work may indeed actively support the liberation of oppressed and marginalized communities, thereby constituting DS4SJ. 
Other work---perhaps the majority---in this area has as its primary goals the furthering of quantitative insight into social systems within the scholarly community. 
An example of the latter kind of work is mathematical analysis of dynamical models of social phenomena, such as segregation or opinion polarization. 
The long-term impact of such work may be primarily localized in scholarly publications, and may not inform activism, policy, or other mechanisms that affect the lives of those who are not professional researchers. 
We emphasize that our aim is not to disparage this body of work, and indeed, several of the present authors maintain active research agendas in this broad area. 

In summary, our focus here---data science for social justice in the mathematical community---overlaps with but is distinct from data science efforts toward broader social good; research on fairness in machine learning; and the quantitative modeling of social systems.

\subsection{The Process of DS4SJ} \label{sec:process}

    Many of us in the mathematical community are accustomed to measuring the success of our research program in units of theorems, experiments, papers, venues, or citations. 
    The processes by which we conduct our research are often optimized towards these measures. 
    Data science for social justice measures its success in terms of impact on the lives of oppressed and marginalized peoples. 
    This focus requires that we pursue the work in ways likely to realize impact, while also minimizing risk of unintended harm.
    The process of this work may therefore look very different from the process of much other mathematical scholarship. 
    We propose a series of four guiding principles for the process of DS4SJ, inspired by \textcite{ardila2016cuentan}.

    \begin{gp}[Collaboration]\label{gp:partnership-collaboration}
        DS4SJ requires partnership with community members and decision-makers in order to realize intended impact and minimize unintended harm. 
    \end{gp}

    \noindent \Cref{gp:partnership-collaboration} reflects two simple facts. 
    First, as mathematicians, we are often inexperienced in conceiving, pursuing, or applying our work outside of academic or industrial research environments.
    We often lack relevant training in social science, policy, or activism. 
    We therefore often need collaboration to execute projects and build skills.
    Second, work that supports a community must be guided by careful engagement \emph{with} that community in order to understand the needs we aim to meet.  
    We identify at least three important elements of effective collaboration for DS4SJ:

        \begin{itemize}
        \item \emph{Perspective} allows us to understand the needs of the people we aim to serve and helps us carry out our work without causing unintended harm.
        \item \emph{Process} helps direct the work with methods that respect the community, and in directions that are responsive to community perspective with results that will support future action.
        \item \emph{Power} makes it possible to actually carry out the intervention informed by our data scientific work. 
    \end{itemize}
    
    Perspective is an especially important element of collaboration, as it helps us understand the needs of the communities we serve.
    The best way to build this perspective is through partnership \emph{with the community itself}. 
    This is ultimately the only reliable way to understand the needs of its members and whether they want what we are able to offer; perspective from other scholars or outside experts can inform only insofar as it is grounded in community engagement. 
    The slogan ``Nothing about us without us,'' popularized in English by South African disability activists, offers a concise maxim.\endnotetext{For more on the origins of the slogan ``Nothing about us without us,'' we refer the reader to \notescitet{charlton2000nothing}}
        
    A recent effort at the Institute for Mathematical and Statistical Innovation (IMSI) illustrates these forms of collaboration. 
    A working group led by author Carrie Diaz Eaton worked in partnership with Nuevas Voces, a leadership program by the Woonasquaucket River Watershed Council (WRWC), for residents of a low-income, immigrant and multilingual community in Providence, Rhode Island.
    This collaboration offered \emph{perspective} by helping the working group understand community needs through the insights of Nuevas Voces, who expressed a need to call attention to the flooding of neighborhoods along the watershed induced by ongoing climate change. 
    The collaborative \emph{process} involved meeting in impacted neighborhoods, conducting interviews in Spanish, and offering support to community members from Nuevas Voces and the WRWC to guide the work as it progressed. 
    Maintaining collaboration through the implementation phase invested the work with \emph{power} to make an impact. 
    The primary deliverable of the intervention with Nuevas Voces was a geospatial dashboard which incorporates not only flooding information but also community resources and stories. 
    Nuevas Voces and WRWC have committed to integrating these tools into their own organizations in order to support their missions.

    \begin{gp}[Critical Reflection]\label{gp:critical-reflection}
        DS4SJ requires critical reflection on one's own identity, position, and privilege, and on the role these play in one's scholarly practice. 
    \end{gp}
    \noindent We use ``critical'' in the sense of \emph{critical theory}: a critical perspective is one which views structures of power and privilege with skepticism, and aims toward a more equitable, free, and just society \autocite{horkheimer1972critical}.
    \endnotetext{For more on critical theory, especially in its connection to liberatory pedagogy, see \notescitet{freire1996pedagogy}.} 
    \Cref{gp:critical-reflection} asks us to reflect on the identities and privileges we inhabit, with special eye towards how these identities and privileges may tacitly inform our work. 
    Many members of the broad mathematical community inhabit multiple forms of professional privilege. 
    Our profession is culturally respected, often commands above-average salaries, and does not usually place us at risk of physical harm.  
    Mathematical ways of knowing are popularly (and often dangerously) viewed as epistemically privileged: piercing, infallible, and not subject to dispute. 
    Many of us in academia have the luxury of measuring deadlines in units of months or years, rather than days or hours. 
    When projects ``don't work out,'' they can often be safely dropped in favor of new ones. 
        Additionally, a disproportionate number of professional mathematicians benefit from White and male privilege.\endnotetext{The disproportionate representation of white men in professional mathematics is highlighted in many places, for example this article in \emph{Scientific American}: \url{https://www.scientificamerican.com/article/modern-mathematics-confronts-its-white-patriarchal-past/}.}
    
    It is important to acknowledge and grapple with privilege when aspiring towards justice work of any form, including DS4SJ. 
    For example, how do we reconcile the typical pace of mathematical research with the fact that, as author Nathan Alexander frequently emphasizes, \emph{people are dying} from systemic forms of injustice?
    The work of interrogating our privileges and understanding how they inform our professional practice is challenging, but there are some aspects of mathematical training that support mathematicians in this work. 
    Privilege can be viewed as a collective social assumption about who deserves freedom, respect, and resources. 
    Assumptions, in the form of theoretical axioms or modeling inputs, are pervasive in mathematical practice. 
    As mathematicians, we are trained in the fundamental skill of articulating assumptions, questioning their applicability, and deriving their consequences. 
    These are exactly the skills that support critical reflection, if they can be turned inward.
    Author Quindel Jones models reflection of the effects of her own identity and positionality on her scholarly work in her vignette in \Cref{sec:vignettes}.

    \begin{gp}[Interrogation of Privilege]\label{gp:divestment-privilege}
        DS4SJ asks the researcher to work actively to interrogate and see past their own privilege. 
    \end{gp}
    As discussed under \Cref{gp:critical-reflection}, mathematicians often benefit from several forms of systemic privilege, often based on a combination of identity, social prestige of profession, and scientific ways of knowing. 
    These forms of privilege can be obstacles to effective DS4SJ work, especially in regards to effective partnership with collaborators from communities that may be very different from our own (\Cref{gp:partnership-collaboration}). 
    Seeing past our our privilege therefore involves epistemic humility: the active centering of ways of knowing other than our own.
    \endnotetext{\notescitet{ho2011trusting} considers epistemic humility in another relationship involving experts, the relationship between physicians and disabled patients. 
    Epistmic humility is also related to \emph{cultural humility}, a concept popularized by \notescitet{tervalon1998cultural}.}
    Learning community needs may require that we heed diverse forms of knowledge, such as community narrative, individual experiences, or common wisdom from other professions.
    Insisting on the primacy of mathematical or scientific ways of knowing over the ways of our served communities is likely to result in failure to synthesize necessary information about community needs or ways to meet them. 
    This can in turn lead to failure to realize intended impact, or even in unintended harm. 
    
    \endnotetext{\notescitet{alexander2021critical} offers further illustrations of interrogating privilege in practice.}
    In his vignette in \Cref{sec:vignettes}, Andr\'es R. Vindas Mel\'endez discusses how his trajectory into DS4SJ as a theoretical mathematician has been aided by awareness of the privilege of academic mathematics and his willingness to engage with other scholarly communities.
                
    \begin{gp}[Systemic Perspective] \label{gp:critical-perspective}
        DS4SJ requires a critical \emph{systemic perspective}: a critical awareness of how systems of power operate to privilege some and oppress others.  
    \end{gp}
    \noindent Many mathematicians use sophisticated tools study the function of complex and interrelated systems in domains ranging from biology to engineering to human society. 
    Systems thinking is a fundamental skill in justice-oriented data work as well, and informs the technical tools we deploy. 
    For example, there is extensive mathematical theory describing mechanisms for fairly distributing resources, but this theory usually lacks a description of historical and present systemic forces of oppression. 
    Rectifying this oppression requires more than fair division of resources; it also requires active work against these systemic forces. 
    In her vignette in \Cref{sec:vignettes}, Heather Zinn Brooks discusses how taking a systemic perspective on the mathematics community can inform data-driven studies of gender representation in mathematical subfields. 
    In some cases, it may also be possible to explicitly model social processes that contribute to disparity. 
    Modeling the social processes of police use-of-force, for example,  may ultimately provide insight on how to stop disproportionate police killings of unarmed Black people \autocite{zhao2020network}.

\section{Challenges and Opportunities in DS4SJ}\label{sec:challenges-opportunities}
\addtoendnotes{\vspace{2em}\noindent {\large \sc Notes on \Cref{sec:challenges-opportunities}}\vspace{1em}}

Justice data science addresses a wide range of problems both within and beyond academia. We roughly divide this work into three overlapping categories. 

\begin{itemize}
    \item \textbf{Data Accessibility}: Ida B. Wells famously wrote that ``the way to right wrongs is to shine the light of truth on them.'' 
    Some justice data science work aims simply to make new kinds of information available and interpretable to new audiences. 
    Examples include efforts to compile accessible, online databases of police violence,\endnotetext{One online database of police violence is the Mapping Police Violence project led by Samuel Sinyangwe: \url{https://mappingpoliceviolence.us}} federal judicial sentencing, 
    \endnotetext{\notescitet{ciocanel2020justfair} obtain a large data set of federal criminal sentencing records, recorded at \url{https://qsideinstitute.org/research/criminal-justice/justfair/}. This data was further analyzed by \notescitet{smith2021racial}.} 
    and cultural representation.\endnotetext{One data set on cultural representation of artists in major US museums was compiled and released by QSIDE at \url{https://qsideinstitute.org/research/the-arts/diversity-of-artists-in-major-u-s-museums-2019/}.}
    \item \textbf{Analysis and Critique}: In some cases, the simple presentation of data in an easily accessible format may be enough to motivate action. 
    In other cases, the sheer quantity of data and the need to account for multiple factors requires explicit, formal data analysis. 
    This analysis can sometimes support controlled or causal claims about the magnitude or impact of bias, and may be published in scholarly journals. 
    In this context, it is important to emphasize that the mere availability of evidence or analysis is often insufficient to change minds or spur change.
    \endnotetext{
        \notescitet{mercier2011humans} offer a theory of reasoning as primarily rhetorical rather than deductive, and includes a survey of work studying confirmation bias and other mechanisms that hinder our ability to change our minds in response to evidence. 
    }
    \item \textbf{Implementation}: Impact is ultimately the axis along which justice work, including data science, is measured. 
    Successful DS4SJ work either achieves impact directly through action or supports follow-up efforts by others.  
    In some cases, it may suffice to bring analysis and critique in front of the right parties; for example, work identifying the racial impact of political gerrymandering is analysis and critique, and can contribute toward its intended impact when presented as arguments in court cases. 
    Other examples include the deployment of predictive or decision-making models aimed at explicit equity goals, such as predicting police misbehavior, landlord exploitation, judicial decision-making, or corporate tax abuse. 
    Data scientists and mathematicians, as engaged human beings, also have opportunities to engage in direct action in support of their chosen causes, including community organizing, attendance at protests, and participation in mutual aid societies.  
    \end{itemize}

\subsection{Getting Involved with DS4SJ}

    We hope that some of our readers are considering whether to engage in or indirectly support work in data science for social justice. 
    We now highlight some of the reasons why this may be appealing for scholars in mathematics and allied fields.  
    There are a wide range of benefits, including active research prospects, ways to engage undergraduates in the classroom, and opportunities to recruit junior collaborators for research projects. 
    
    Justice data science work can motivate the development of new 
    mathematics, data scientific techniques, and methodological perspectives.
        One well-known example is the development of methods for detecting and quantifying racial bias in electoral maps by defining null distributions on spaces of districting schemes using techniques from metric geometry and using Markov-chain algorithms to sample from them  \autocite{becker2021computational}. 
    This body of work has produced analysis and critique in the form of papers quantifying racial bias in existing political maps. 
    There have also been impacts from engagement with power partnerships, with researchers collaborating with state governments\endnotetext{Moon Duchin offers a popular description of collaborations between metric geometers, legislators, and legal experts in a recent article in \emph{Scientific American}: \url{https://www.scientificamerican.com/article/geometry-versus-gerrymandering/}} and contributing to amicus briefs before the US Supreme Court.\endnotetext{One amicus brief submitted by the Metric Geometry and Gerrymandering Group to the Supreme court may be found at \url{https://mggg.org/SCOTUS-MathBrief.pdf}}  

    There are also opportunities for methodological advancements in applied statistics, optimization, and network science.
    \textcite{pierson2020large}, for example, develop novel, highly-scalable statistical tests for detecting racial bias in large-scale datasets of police behavior.\endnotetext{\notescitet{pierson2018fast} offer additional work on large-scale statistical tests for bias.}
    Problems inspired by equitable access to schools, grocery stores, and polling places in urban infrastructure networks have led to the development of novel algorithms for editing graph edges to influence diffusion dynamics \autocite{ramachandran2021gaea}.
    A growing body of work uses network-based and statistical methods to quantify gendered inequalities in citation, promotion, and retention across academic disciplines \autocite{wapman2022quantifying}.
    In her vignette below, author Heather Zinn Brooks describes ongoing related work in mathematics.
    As author Andr\'es R. Vindas Mel\'endez writes in his own vignette, the urgency of social justice applications have helped him both develop an interest in data science techniques and begin to reflect on whether combinatorics, his primary area of research training, might have useful applications toward justice causes. 
    
    An orientation towards social justice offers benefits for pedagogy in mathematics, data science, statistics, and computer science.
    Real-world topics related to issues of injustice and inequity are fruitful sources of problems and examples \autocite{karaali2019mathematics}. 
    Beyond simply furnishing problems, however, justice orientations can support mathematics instructors in welcoming, engaging, and retaining students in undergraduate STEM majors/minors, and higher education more generally. 
    Justice orientations can be especially useful for engaging students who have been historically excluded from or oppressed or marginalized within traditional classrooms. 
    \endnotetext{\notescitet{garibay2015stem} discusses justice orientation as a promoter of student engagement.}
    Connecting classroom content to ``prosocial communal outcomes'' is a recommendation of a recent report focused on improving the persistence of marginalized students in STEM disciplines \autocite{estrada2016improving}. 
    
    It is important to note that justice-oriented pedagogy presents challenges for both students and instructors, and can be harmful to students if cases are presented in a demeaning or tokenizing way. 
    There is potential for harm if instructors incorporate examples that do not promote positive student identity or which fail to acknowledge diversity of experience within cultures and identities \autocite{leonard2010nuances}.

    Effective justice-oriented pedagogy seeks to advance students' critical consciousness--their ability to challenge what is framed as normal, whose interests are served by that framing, how system structures are reified, and in what ways they could, potentially, be different.  
    Such pedagogy not only welcomes students who have been historically marginalized in classroom spaces, but also challenges all students to grow as justice-oriented systems thinkers. 
    \endnotetext{For further suggestions on critical conversations to support justice-oriented mathematics pedagogy, see \notescitet{alexander2021critical}.}
    For example, author Nathan Alexander uses a justice orientation in his introductory statistics classroom when framing the topic of ethical research and data use. 
    The Institutional Review Board (IRB) system was formed in part due to the revelation of the medical abuse of almost 400 Black men in the Tuskegee Syphilis Study.
    The modern IRB review system aims to prevent future such flagrant breaches of medical ethics. 
    This, however, does not fully address the dimensions of structural racism in the practice of medical research; even today, Black principal investigators are less likely to have grant proposals funded by the National Institutes of Health. \endnotetext{Disparity in funding of NIH proposals is studied by \notescitet{taffe2021equity}.} 
    These considerations are not only context for the modern practice of statistics and data-informed research, but they are also opportunities to practice core technical course content. 
    Students can study and evaluate, for example, the methods underlying recent findings that the Tuskegee study is \emph{still} shortening the lives of Black men by sowing justified mistrust in medical institutions \autocite{alsan2018tuskegee}. 

    These are challenging topics to navigate without tokenizing or otherwise harming students. 
    Doing so requires instructors to show vulnerability, humility, and authentic willingness to learn \emph{with} students in a community built on mutual respect, openness, and collective critical reflection. 
    The reward is an opportunity for students to develop technical skills and critical awareness in tandem, each reinforcing the other. 
    \endnotetext{For further critical perspective on the role of justice and equity-oriented data science, see \notescitet{alexander2022beyond} and \notescitet{castillo2022}.
    \notescitet{eaton2022teaching} offers a detailed description of implementing justice orientation in machine learning pedagogy.}

    Similarly, students at all levels are often excited to work on research projects that involve justice work. 
    Undergraduates are increasingly involved in and committed to social justice and are drawn to research experiences with obvious real-world applications \autocite{seemiller2017generation}. 
    As author Ariana Mendible describes in her vignette below, her early research projects with justice focus have received much higher than expected interest from prospective research students at her institution. 
    This may be in part due to the accessibility of these projects. 
    DS4SJ work often has lower technical barrier to entry when compared to many other problems in academic mathematics, making it possible to involve students across a greater range of experiences and academic backgrounds at earlier stages in their academic careers. 
    Students' personal experiences and interests can also be explicit assets in the pursuit of DS4SJ research projects, especially when they are members of the communities who stand to be impacted by the work.

        Finally, we emphasize that DS4SJ work places distinctive demands on researchers and instructors, and is therefore not necessarily for everyone. 
        First, scholars may reasonably have motivations and interests that are not fully aligned with DS4SJ work. 
        DS4SJ work measures itself first and foremost on positive impact for marginalized or oppressed peoples. 
        Scholars pursuing this work therefore need both orientation towards concrete outcomes and willingness to devote time towards realizing those outcomes. 
        This in turn requires time to cultivate critical perspective, engage carefully with communities, and develop nontechnical skills. 
        We have nothing but respect for scholars who reflect on the demands of DS4SJ and decide to pursue other paths. 
        Second, we acknowledge that professional incentives in many institutions and societies may not encourage DS4SJ work. 
        For example, some departments may not recognize DS4SJ scholarship as being ``sufficiently mathematical'' to warrant promotion or recognition, or may view justice-oriented content in the classroom as a distraction. 
        Scholars cannot be expected to commit to this work if they feel it endangers their prospects for advancement in their departments, institutions, or communities.

\subsection{DS4SJ and the Math Community}

    We therefore call on departments, institutions, and professional communities to support work in DS4SJ. 
    Fortunately, DS4SJ work is aligned with the stated goals of these departments and communities, implying both opportunity and incentive to support these efforts. 
    The AMS Statement on Equity, Diversity, and Inclusion states our shared goal of \emph{a community that is diverse, respectful, accessible, and inclusive}. 
    We are still far from this goal in 2023. 
    Academic mathematics still has not reached proportional representation or equity of opportunity along the lines such as gender and race.\endnotetext{\url{https://www.scientificamerican.com/article/modern-mathematics-confronts-its-white-patriarchal-past/}}
    Many departments struggle to recruit and retain minoritized trainees and faculty. 
    Active support of justice data science may be able to contribute to these goals. 
    Based on interviews conducted with underrepresented racial minority students in US STEM graduate programs, \textcite{tran2011can} recommends that departments seeking to retain underrepresented minority students should legitimize and expand research that has direct applications for social justice. 
    Illustrating this effect, author and Ph.D. candidate Quindel Jones  discusses how the opportunity to use mathematical models to improve health outcomes for Black women is a motivating thread in her career trajectory in her vignette below. 
    Jones highlights how the presence of a growing community of mathematical scholars focused on data science and social justice has promoted her scholarship, enabled her to learn from mentors who shared some of her identities, and connect with collaborators who share her interests. 
    Support of justice data science by departments, research institutions, and professional organizations can not only promote scholarship that stands on its own merits, but also further our shared mission towards a more equitable, diverse, and inclusive mathematics profession. 
    We discuss concrete ways that departments,  institutes, and societies can support DS4SJ in \Cref{sec:conclusion}.

\section{Experiences in DS4SJ} \label{sec:vignettes}

\addtoendnotes{\vspace{2em}\noindent {\large \sc Notes on \Cref{sec:vignettes}}\vspace{1em}}

    In this section, we highlight the experiences of several pre-tenure mathematicians in their diverse engagements with data science for social justice. 
    Their stories offer examples of how the principles of DS4SJ described in \Cref{sec:process} can be implemented, and how the benefits and challenges of this work have influenced their trajectories. 
    Ariana Mendible discusses how collaboration has impacted and improved her DS4SJ work (\Cref{gp:partnership-collaboration}). 
    Quindel Jones provides a critical reflection of identity and positionality on scholarship (\Cref{gp:critical-reflection}). 
    Andr\'es R. Vindas Mel\'endez interrogates the privilege and challenges of doing DS4SJ work as a theoretical mathematician (\Cref{gp:divestment-privilege}). 
    Heather Zinn Brooks highlights the importance of taking a systemic perspective on research work in DS4SJ (\Cref{gp:critical-perspective}). 
    Their experiences highlight the roles of identity, belonging, and justice orientation in the trajectories of early-career scholars working in science for social justice. 

\subsubsection*{Ariana Mendible, Seattle University}

    I recently returned to my alma mater, an institution with a strong commitment to justice, as an assistant professor of mathematics.
    I recognize the privilege of this position, offering the freedom to work on unfunded problems that align with my personal values.
    Pivoting my research to DS4SJ has offered fruitful research collaborations that promise meaningful impact within and outside of academia.

    My work as a co-director of QSIDE's Small Town Police Accountability (SToPA) research lab has shown me the eagerness of rising scholars to engage with DS4SJ. 
        We have worked together to tackle tasks like optical character recognition and topic modeling with the goal of empowering small towns to access and understand data from their police departments. The impact and urgency of these problems invites a broad audience to participate in research, including students who may not otherwise have envisioned themselves as researchers.

    Sociologists, lawyers, and activists participating in our lab efforts enrich our perspectives beyond our mathematical and technical expertise. 
        Importantly, collaboration with community leaders guarantees direct and effective change from our work. Conversations with one town's racial justice and police reform group provided town-specific history and helped narrow our focus to the most impactful questions. 
        Community members have expressed enthusiasm that our quantitative work will further empower them to reform their local police departments. 
    We hope that their enthusiasm will help our work inform community decision-making and realize its intended impact beyond the ivory tower.

    My long-term goal is to continue building a research program that invites students to work on problems that are meaningful to them and their communities. I believe that this program will not only enact my values, but ultimately strengthen my academic career.

\subsubsection*{Quindel Jones, Virginia Commonwealth University}

    I am a queer Black woman born and raised in Jackson, Mississippi, a part of the American South. 

    While I initially took to dance and writing, my mother's affinity to numbers engaged my mathematical curiosity. 
    My own affinity grew as community summer enrichment programs kept me engaged. 
    Throughout high school and college, there were opportunities for me to use mathematics in my personal life and learn from people who looked like or cared about me. 
    I graduated from Jackson State University, one of six historically Black colleges and universities (HBCUs) in the state, under the guidance of Dr. Jana Talley, whom I met in the I.C.STEM program at Jackson State when I was in high school. 

    I am currently a doctoral candidate in applied mathematics with a concentration in mathematical biology at Virginia Commonwealth University. 
    Since entering graduate school, my professional mission has been to use mathematics to improve the lives of people that look like me---especially Black women. 
    In my dissertation project, I am developing a dynamical model of pain levels during sickle cell disease, a disease that disproportionately affects Black people. 
    We start with self-report data of patient pain, collected by a member of our interdisciplinary team. 
    We use this data to fit an ordinary differential equation (ODE) model of a patient's pain profile over time. 
    This model is already informing the data collection strategies for the next round of our study. 
    Long-term, we hope to embed our model in a wearable app. 
    We hope for patients to be able to log their sleep and pain self-reports, and from these receive warnings about likely upcoming pain crises. 
    Patients could use these warnings to seek prophylactic care or medication. 

    As I continue on my mathematics journey, I am discovering more and more ways in which mathematics can be applied to address inequities. 
    With my own life having been shaped by inequity due to my marginalized identity in U.S. society, I am thrilled that there's an active community centered on math and data science for social justice. 
    Connecting with this community in venues like the recent program at ICERM have been formative experiences in my mathematical career. 
    These programs and their leaders have helped me see myself and my concerns as important. 
    Part of my excitement for my future as a researcher is knowing that socially impactful mathematical work exists, and that there is a growing community of collaborators with whom to do it. 

\subsubsection*{Andr\'es R. Vindas Mel\'endez, UC Berkeley}

    My mathematical work is mainly in  combinatorics. 
    I have always found value in this work from meaningful interactions with collaborators and the joy of mathematical discovery.

    I had my ``formal'' introduction to data science for social justice at a recent ICERM program on this topic. 
    I began working with a number of collaborators to address a deceptively simple question: \emph{who} or \emph{what} is a mathematician?
    Narrow definitions based on job titles or institutional affiliations can be precise, but can also be used to implicitly or explicitly exclude people from the discipline. 
    In a recent preprint, we examined some of the tensions in several formulations of what it means to be a mathematician, including those based on self-identification, activities, and qualifications.
    \endnotetext{The preprint referred to by Andr\'es R. Vindas Mel\'endez is: \notescitet{defs}}
    For future work, we aim to design a survey of mathematicians (broadly construed) across a wide range of institutions, personal identities, job titles, and areas of scholarly and pedagogical interest. 
    We will study these data using techniques like text analysis and  clustering to extract qualitative insight from survey participant responses. 
    From this work, we hope to deepen our understanding of the identities and needs of mathematicians in order to support interventions at multiple levels to increase the participation of marginalized or underrepresented groups. 

    My interest in the topic of who is considered a mathematician, and who is welcomed in our community, is informed by my own identity and experience. 
    My identities as a queer, chronically-ill, first-generation,
        Latinx mathematician have shaped my own experience in academic mathematics, including both successes and struggles.
    I now benefit from privilege in my study and work and well-resourced institutions. 
    My experiences navigating my career so far highlight the need for spaces that allow marginalized people to feel welcomed, comfortable, and supported in the mathematics community. 

    As someone whose training is primarily in theoretical math, I sometimes feel that I do not have much to contribute to DS4SJ. 
    I remind myself that this work requires diverse experience and expertise and my own identities offer important perspective on critical questions.
    What is missing from our survey efforts? 
    What identities might be overlooked or mischaracterized by our approach? What voices are centered in extant literature? 
    Whose voices are neglected?  
    
    Recently, I have been developing my technical data science skills. 
    At my current institution, 
        I am participating in workshops on Python, R, data visualization, and text analysis as a means to build my skills for my research projects. 
    Furthermore, I have joined the Digital Humanities Working Group, which allows me to have interdisciplinary conversations and reflect on the social prestige that mathematicians carry.

\subsubsection*{Heather Zinn Brooks, Harvey Mudd College}
    
    My identities as a woman and a first-generation college student have impacted my career trajectory 
    and inspired my involvement in research in data science and social justice.
    Growing up, I didn't know that ``mathematician'' was a real job title.
    As a student, I struggled to see myself as a mathematician because of the lack of role models who shared these identities.
    Now that I have the privilege to carry the title of ``professor of mathematics,''
    one of my motivations for becoming involved in work in DS4SJ is to give back to (and improve) a community that has benefited me. I recognize that my positionality as a White, cis-gendered, able-bodied individual has contributed toward my positive experiences in the mathematics community. 
    My interest in creating inclusive mathematical spaces has led me to a recent project modeling female gender representation in mathematical subfields and institutions.
    In this project, my collaborators and I use a data set compiled from the Mathematics Genealogy Project (MGP), combined with algorithmically inferred binary gender from names. 
    This data has important limitations. 
        Binary categories are imperfect representations of gender, an inherently nonbinary facet of identity.
        Furthermore, gender is not determined by name or appearance. 
        We assume that misgendering takes place in our data set, and we view our results as at most estimates of macroscopic trends.
    \endnotetext{Another important limitation of the automated gender inference methods described in Heather Zinn Brooks' vignette is the Euro- and US-centrism of most commercially available algorithms, including the ones we used. 
        For a review of these biases and limitations, see \notescitet{santamaria2018comparison}} 
    Within these limitations, our goal is to build understanding of how some subfields and institutions have reached relatively high levels of female gender representation, while others remain strongly male-dominated. 

            We study the roles of accumulation, homophily, and prestige as mechanisms supporting or inhibiting female gender representation in mathematics. 
    A centerpiece of our work is the development of data-informed branching process models for replicating features found in real data. 
    We use machine learning techniques to highlight useful predictors in the data and inform the model mechanisms. 
    In addition to these technical contributions, a long-term aim is to develop actionable strategies to support more inclusive mathematical spaces. 
    
     My formal training in applied dynamical systems and mathematical biology has allowed me to connect this motivation into my research in meaningful ways and has challenged me to adopt a systemic perspective on the way our mathematical community functions. 
          Personal identity characteristics, institutional prestige, and perceptions of and the microcultures within mathematical subfields are all salient features that impact the way our community looks, feels, and operates. This systemic perspective allows my collaborators and I to interrogate the MGP data to highlight, explore, and understand the impact of these features.
     
    My primary motivation for engaging in justice-oriented data science work is to fight systemic injustice and oppression.
    My involvement in this work has also inspired new scholarship, provided a source of joy in my research career, and allowed me to connect my identity to my mathematical pursuits.    

\section{Outlook} \label{sec:conclusion}

\addtoendnotes{\vspace{2em}\noindent {\large \sc Notes on \Cref{sec:conclusion}}\vspace{1em}}

    We have argued that many members of the mathematical community can, if they choose, further justice work through the reflective, critical practice of data science. 
    We have highlighted some of the benefits for scholars, especially early-career mathematicians: 
    DS4SJ allows scholars to do work with concrete impact, to relate their identities to their work, to find inspiration for new technical problems, and to connect with their students inside and outside the classroom. 
    There are communities emerging to support this growing subdiscipline, foster collaborations, and build critical expertise.  
    We close with some suggestions for individuals, departments, and institutions to nurture work in this burgeoning area.

While individual motivation may spark engagement in DS4JS scholarship, departments, institutions, and professional societies must offer active support for this scholarship to be successful. 
This is especially pertinent considering that a large portion of the work in DS4SJ is being spearheaded by early-career mathematicians, including untenured professors, industry professionals, postdocs, and students. 
Since many of these scholars lack career stability, it is important that institutions align incentives so that these scholars can be confident in their ability to advance their careers while pursuing justice work. 
At all levels, DS4JS work must be structurally and financially supported to sustain scholarship in this field. 

Academic institutions can support justice-oriented scholars at all career stages.  
Departments can support faculty in their teaching endeavors by encouraging the development of new courses at the intersection of mathematics, data science, and social justice. 
Programs can provide training to undergraduate and graduate students and facilitate opportunities for community collaborations. 
Hiring committees can revise rubrics to explicitly place value on justice-oriented work, which may be particularly important in searches related to applied data science. 
In retention, promotion, and tenure decisions, work must be understood and valued appropriately. 
For many departments, this necessitates broadening the perspective on what constitutes evidence of scholarship in this area, including both publications and evidence of impact of the scholar’s work outside academic venues. 
In order to meaningfully support justice-focused scientific work, departments must formally recognize its value in training and evaluative processes at all levels.

The structural support of institutes and professional societies is also critical to furthering this growing subfield of mathematics. 
As with any subfield, workshops and conferences are critical for disseminating work, generating new ideas, and creating and maintaining connections among scholars. 
Such events can be especially impactful for early-career scholars aiming to make connections and jump-start their scholarship. 
The Summer 2022 program on ``Data Science and Social Justice: Networks, Policy, and Education'' at ICERM, at which the present authors began this article, is just one example of a growing number of research community spaces.

\endnotetext{Other examples include a recent workshop on Mathematics and Racial Justice at the Mathematical Sciences Research Institute (MSRI), the upcoming workshop on Interdisciplinary and Critical Data Science Motivated by Social Justice at the Institute for Mathematical and Statistical Innovation, and the second part of the program on Data Science and Social Justice at the Institute for Computational and Experimental Research in Mathematics. We also note the long-standing annual Critical Inquiry in Mathematics Education conference at MSRI.}

Simply making venues available for workshops and collaborations is, however, not sufficient. 
Adequate, equitable funding is necessary for diverse researchers and community members to fully participate in partnership. 
Providing this funding may challenge traditional institutions. 
Some examples of funding to support equitable participation include travel for community participants (not just researchers); upfront financial support to graduate student researchers for travel and meals; and financial support for childcare or partner travel. 
Measures like these enable the full participation of individuals for whom the typical financial and time costs of workshop participation may be especially heavy.
We therefore encourage reflection and creativity on the part of institutes and other funded projects to design funding which supports equitable participation in DS4SJ work. 

Professional societies can further research and innovation by forming and supporting professional communities for mathematicians working in DS4SJ. 
Messaging from these societies is a powerful and important tool for building respect and professional legitimacy for justice data scientists and the scholarship they create. 
Such messaging also helps scholarly reviewers understand the nature of DS4SJ research, which this article aims to inform.

The intersection of data science and social justice work represents an exciting opportunity for mathematical scientists to work toward a more just and equitable future, both within our own community and beyond. 
We hope that this article can promote the growth and success of data science for social justice. 

\section*{Author Positionality}

    The authors of this manuscript represent a diverse range of identities, experiences, and voices within the mathematical community. 
    Our views on data science and social justice are necessarily shaped by our experiences of both marginalization and privilege. 
    Many of us embody marginalized identities along axes of race, ethnicity, gender, and sexuality. 
    For some of us, our interest in justice work is informed by a passion to combat forms of oppression that we ourselves have experienced. 
    We also acknowledge that many of us occupy positions of privilege within the mathematical community, including job security, access to research funding, and recognition within our professional communities. 
    We look forward to continuing both our learning and our work towards justice.

\section*{Acknowledgements}

    This material is based upon work supported by the National Science Foundation under Grant No. DMS-1929284 while some the authors were in residence at the Institute for Computational and Experimental Research in Mathematics in Providence, RI, during the Summer 2022 program on ``Data Science and Social Justice: Networks, Policy, and Education.'' 
    Part of this research was performed while some of the authors were visiting the Institute for Mathematical and Statistical Innovation (IMSI), which is supported by the National Science Foundation (Grant No. DMS-1929348).
    Vindas-Mel\'endez is partially supported by the National Science Foundation under Award DMS-2102921.

\printbibliography[notcategory=fullcited]

\onecolumn

\setcounter{page}{1}

\newpage 

\theendnotes

\printbibliography[category=fullcited,title={References in Notes}]

\end{document}